\documentclass[12pt,draft]{article}

\usepackage[cp1251]{inputenc}
\usepackage[T2A]{fontenc}
\usepackage[english,russian]{babel}
\usepackage{amsmath,amsfonts,amsthm,amssymb}
\usepackage{geometry}
\theoremstyle{plain}

\theoremstyle{definition}

\usepackage[space]{cite}

\begin{document}

\begin{flushleft}

\medskip

\textbf{A. A. Murach }\small{(Institute of Mathematics of the National Academy of Sciences of Ukraine, Kyiv)}
\normalsize

\textbf{O. B. Pelekhata }\small{(National Technical University of Ukraine \glqq Igor Sikorsky Kyiv Polytechnic\\ \hspace{3.5cm} Institute\grqq, Kyiv)} \normalsize

\textbf{V. O. Soldatov }\small{(Institute of Mathematics of the National Academy of Sciences of Ukraine, Kyiv)} \normalsize

\medskip

\large\textbf{Approximation properties of solutions\\to multipoint boundary-value problems}

\end{flushleft}

\begin{flushleft}

\bigskip

\textbf{А. А. Мурач }\small{(Институт математики НАН Украины, Киев)}
\normalsize

\textbf{О. Б. Пелехата }\small{(Национальный технический университет Украины \glqq Киевский \\ \hspace{3.3cm} политехнический институт имени Игоря Сикорского\grqq, Киев)} \normalsize

\textbf{В. А. Солдатов }\small{(Институт математики НАН Украины, Киев)} \normalsize

\medskip

\large\textbf{Аппроксимационные свойства решений\\многоточечных краевых задач}

\end{flushleft}

\bigskip

\noindent We consider a wide class of linear boundary-value problems for systems of $m$ ordinary differential equations of order $r$, known as general boundary-value problems. Their solutions $y:[a,b]\to \mathbb{C}^{m}$ belong to the Sobolev space $(W_1^{r})^m$, and the boundary conditions are given in the form $By=q$ where $B:(C^{(r-1)})^{m}\to\mathbb{C}^{rm}$ is an arbitrary continuous linear operator. We prove that a solution to such a problem can be approximated with an arbitrary precision in $(W_1^{r})^m$ by solutions to multipoint boundary-value problems with the same right-hand sides. These multipoint problems are built explicitly and do not depend on the right-hand sides of the general boundary-value problem. For these problems, we obtain estimates of errors of solutions in the normed spaces $(W_1^{r})^m$ and $(C^{(r-1)})^{m}$.

\vspace{0,5cm}

\noindent\textbf{1. Введение.} Предельный переход в системах дифференциальных уравнений часто используется при решении различных математических и прикладных задач, поэтому его обоснованию и исследованию его свойств посвящено много работ. Наиболее полно эта проблематика изучена  для задачи Коши для дифференциальных уравнений первого порядка. Значительно менее исследованы краевые задачи, зависящие от параметра. Такая ситуация связана с большим разнообразием краевых условий. Первые результаты в этом направлении получили И.~Т.~Кигурадзе \cite{Kiguradze1987, Kiguradze1975, Kiguradze2003} и М.~Ашордия \cite{Ashordia1996} для широкого класса так называемых общих линейных краевых задач для систем обыкновенных дифференциальных уравнений первого порядка. Были найдены достаточные условия непрерывности по параметру решений этих задач в нормированном пространстве $C([a,b],\mathbb R^m)$, где $m$~--- количество уравнений в системе. В работах В.~А.~Михайлеца и его учеников \cite{KodliukMikhailetsReva2013, MikhailetsReva2008DAN9, MikhailetsChekhanova2015JMathSci, MikhailetsPelekhataReva2018UMJ, PelekhataReva2019UMJ} эти результаты были существенно улучшены и обобщены на случай комплекснозначных функций и систем дифференциальных уравнений высокого порядка~$r$. Решения этих систем пробегают пространство Соболева $(W_1^{r})^m=W_1^{r}([a,b],\mathbb C^{m})$, а краевые условия имеют вид $By=q$, где $B$~--- произвольный непрерывный линейный оператор на паре пространств $(C^{(r-1)})^m=C^{(r-1)}([a,b],\mathbb C^{m})$ и $\mathbb C^{rm}$. Другие широкие классы линейных одномерных краевых задач введены и исследованы в работах \cite{MikhailetsReva2008DAN8, KodliukMikhailets2013JMS, Soldatov2015UMJ, Gnyp2016, MikhailetsMurachSoldatov2016EJQTDE, MikhailetsMurachSoldatov2016MFAT, GnypKodlyukMikhailets2015UMJ, GnypMikhailetsMurach2017}. Решения этих задач рассматриваются в нормированных пространствах $n$ раз непрерывно дифференцируемых функций, где $n\geq r$, функциональных пространствах Гельдера порядка $l>r$ или Соболева\,--\,Слободецкого порядка $l\geq r$, причем на этих же пространствах заданы и непрерывные краевые операторы.

Важное место среди упомянутых задач занимают многоточечные краевые задачи. Поскольку краевые условия в этих задачах имеют достаточно простой явный вид, их решения найти гораздо легче, чем, например, решения задач с интегральными краевыми условиями. Поэтому возникает естественный вопрос о возможности аппроксимации решениями многоточечных краевых задач решений более общих краевых задач. Цель настоящей работы~--- показать, что решения указанных выше общих краевых задач можно с произвольной точностью аппроксимировать в пространстве $(W_1^{r})^m$ решениями многоточечных краевых задач с теми же правыми частями. Эти аппроксимирующие задачи строятся нами явно и не зависят от правых частей общей краевой задачи. Для этих задач мы получаем оценки погрешности решений в нормированных пространствах $(W_1^{r})^m$ и $(C^{(r-1)})^m$. Для другого класса линейных краевых задач, решения которых принадлежат пространству $(C^{(n)})^{m}$ целого порядка $n\geq r$, соответствующий результат об аппроксимации в $(C^{(n)})^{m}$ доказан в недавней статье \cite{MasliukPelekhataSoldatov2020MFAT}.

\textbf{2. Постановка задачи и результаты.} Пусть $a,b \in\mathbb{R}$ и $a<b$. Используем такие краткие обозначения комплексных банаховых пространств функций $x:[a,b]\to\mathbb{C}$ и их норм:
\begin{itemize}
  \item[(a)] $L_1$~--- пространство всех интегрируемых по Лебегу на $[a,b]$ функций с нормой $\|x\|_1:=\int_a^b|x(t)| dt$;
 \item[(b)] $C$ или $C^{(0)}$~--- пространство всех непрерывных на $[a,b]$ функций с нормой $\|x\|_{C}:=\|x\|_{(0)}:=\max\limits_{a\leq t\leq b}|x(t)|$;
 \item[(c)] $C^{(l)}$, где $l\in\mathbb{N}$,~--- пространство всех $l$ раз непрерывно дифференцируемых на $[a,b]$ функций с нормой $\|x\|_{(l)}:=\sum_{j=0}^{l}\|x^{(j)}\|_{C}$;
  \item[(d)] $W_1^r$, где $r\in\mathbb{N}$,~--- пространство Соболева всех функций $x\in C^{(r-1)}$ таких, что функция $x^{(r-1)}$ абсолютно непрерывна на $[a,b]$ (тогда $x^{(r)}(t)$ существует для почти всех $t\in[a,b]$, причем $x^{(r)}\in L_1$); в пространстве $W_1^r$ введена норма $\|x\|_{r,1}=\sum_{j=0}^{r}\|x^{(j)}\|_1$.
\end{itemize}

Комплексные банаховы пространства, образованные вектор-функциями размерности $m\geq\nobreak1$ или матрицами-функциями размера $m\times m$, все компоненты которых принадлежат одному из перечисленных пространств, обозначаем соответственно через $(\Psi)^{m}$ или $(\Psi)^{m\times m}$, где $\Psi$ символизирует одно из указанных пространств скалярных функций. При этом векторы интерпретируем как столбцы. Норма вектор-функции в пространстве $(\Psi)^{m}$ равна сумме норм всех ее компонент в $\Psi$, а норма матрицы-функции в пространстве $(\Psi)^{m\times m}$ равна максимуму норм всех ее столбцов в $(\Psi)^{m}$. Нормы в пространствах $(\Psi)^{m}$ и  $(\Psi)^{m\times m}$ обозначаем таким же образом, как и норму в пространстве $\Psi$. Из контекста всегда будет понятно о норме в каком пространстве (скалярных функции, вектор-функций или матриц-функций) идет речь. Для числовых векторов и матриц используем аналогичные нормы, которые  обозначаем через $\|\cdot\|$.

Пусть $r,m\in\mathbb{N}$. На отрезке $[a,b]$ рассмотрим линейную краевую задачу вида
\begin{align}\label{72MSp.syste}
Ly(t):=y^{(r)}(t)+\sum_{l=1}^r A_{r-l}(t)\,y^{(r-l)}(t)&=f(t)\quad \mbox{для п.в.}\;\;t\in[a,b],\\
By&=q.\label{72MSp.kue}
\end{align}
Здесь вектор-функция $y$ класса $(W_1^{r})^m$ искомая, и произвольно заданы матрицы-функции $A_{r-l}\in(L_1)^{m\times m}$, где $1\leq l\leq r$, вектор-функция $f\in(L_1)^{m}$, вектор $q\in\mathbb{C}^{rm}$ и непрерывный линейный оператор
\begin{equation}\label{72MSp.oper}
B:(C^{(r-1)})^{m}\rightarrow \mathbb C^{rm}.
\end{equation}
Из условия $y\in(W_1^{r})^m$ следует, что старшая производная $y^{(r)}$ в уравнении \eqref{72MSp.syste} существует почти везде на $[a,b]$ (относительно меры Лебега) и является вектор-функцией класса $(L_1)^{m}$. Поэтому это уравнение рассматривается для почти всех (п.в.) $t\in[a,b]$. Поскольку пространство $(W_1^{r})^m$ непрерывно вложено в $(C^{(r-1)})^{m}$, краевое условие \eqref{72MSp.kue} имеет смысл. Такую задачу исследовали в \cite{MikhailetsChekhanova2015JMathSci, MikhailetsPelekhataReva2018UMJ, PelekhataReva2019UMJ}.

Далее предполагаем, что краевая задача \eqref{72MSp.syste}, \eqref{72MSp.kue} имеет единственное решение $y\in(W_1^{r})^m$ для произвольных правых частей $f\in(L_1 )^{m}$ и $q\in\mathbb C^{rm}$. Согласно \cite[с.~334]{MikhailetsChekhanova2015JMathSci}, это предположение эквивалентно тому, что соответствующая однородная краевая задача (для $f(\cdot)\equiv0$ и $q=0$) имеет лишь тривиальное решение $y(\cdot)\equiv0$.

Пусть $\mathcal{X}$~--- произвольное плотное подмножество пространства $(L_1)^{m\times m}$ (например, множество всех полиномов или тригонометрических полиномов, ступенчатых функций, полигональных функций,  сплайнов). Рассмотрим последовательность многоточечных краевых задач вида
\begin{align}\label{73MSp.syste}
L_k\,y_k(t):=y_k^{(r)}(t)+\sum_{l=1}^r A_{r-l,k}(t)\,y^{(r-l)}_{k}(t)&=f(t)\quad\mbox{для п.в.}\;\;t\in[a,b],\\
B_{k}\,y_k:=\sum\limits_{j=1}^{p_k}\sum
\limits_{l=0}^{r-1}\beta_k^{j,l}\,y^{(l)}(t_{k,j})&=q. \label{7MSp.2B_m}
\end{align}
Они параметризованы натуральным числом $k$; их правые части не зависят от $k$ и совпадают с правыми частями краевой задачи \eqref{72MSp.syste}, \eqref{72MSp.kue}. Относительно этих многоточечных задач предполагаем, что каждое $A_{r-l,k}\in \mathcal{X}$ и, кроме того, $p_k\in\mathbb{N}$, $\beta_k^{j,l}\in \mathbb{C}^{rm\times m}$ и $t_{k,j} \in [a,b]$  для всех допустимых значений параметров $k$, $j$ и~ $l$. Решения $y_k$ этих задач рассматриваются в классе $(W_1^{r})^m$. Очевидно, линейный оператор $B_k$ действует непрерывно из $(C^{(r-1)})^{m}$ в $\mathbb{C}^{rm}$.

Сформулируем основной результат работы.

\textbf{Теорема 1.} \it Для краевой задачи \eqref{72MSp.syste}, \eqref{72MSp.kue} существует последовательность многоточечных краевых задач вида \eqref{73MSp.syste}, \eqref{7MSp.2B_m} таких, что они однозначно разрешимыми для достаточно больших $k$ и выполняется асимптотическое свойство
\begin{equation}\label{7MSp.2granum}
y_{k}\to y \;\;\mbox{в}\;\; (W^r_1)^m\;\;\mbox{при}\;\;k\rightarrow\infty.
\end{equation}
Эту последовательность можно выбрать независимой от $f$ и $q$ и построить явно. \rm

Далее в этом разделе предполагаем, что последовательность краевых задач \eqref{73MSp.syste}, \eqref{7MSp.2B_m} удовлетворяет заключению теоремы~1 для произвольных правых частей $f\in(L_1)^{m}$ и $q\in\nobreak\mathbb{C}^{rm}$, т.~е. эти задачи однозначно разрешимы в пространстве $(W^r_1)^m$ при $k\geq\widetilde{\varrho}$ и их решения имеют свойство~\eqref{7MSp.2granum}; здесь $\widetilde{\varrho}$~--- некоторое натуральное число.

Исследуем случай, когда правые части этих задач зависят от $k$. Таким образом, при $k\geq\widetilde{\varrho}$ рассмотрим краевые задачи вида
\begin{align}\label{Lk=fk}
L_{k}\,x_k(t)&=f_k(t)\quad\mbox{для п.в.}\;\;t\in[a,b],\\
B_{k}\,x_k&=q_k, \label{Bk=gk}
\end{align}
где $f_k\in(L_1)^{m}$ и $q_k\in\mathbb{C}^{rm}$, а решение $x_k\in(W^r_1)^m$.

\textbf{Теорема 2.} \it Пусть заданы натуральное число $\widehat{\varrho}\geq\widetilde{\varrho}$ и действительное число $\varepsilon>0$. Предположим, что правые части краевых задач \eqref{Lk=fk}, \eqref{Bk=gk}
удовлетворяют условиям
\begin{equation}\label{74aut.r.s.est}
\|f_k-f\|_{1}<\varepsilon\quad\mbox{і}\quad\|q_k-q\|<\varepsilon
\quad\mbox{при}\;\;k\geq\widehat{\varrho}.
\end{equation}
Тогда существуют положительные числа $\varkappa$ и $\varrho\geq\widehat{\varrho}$ такие, что
\begin{equation}\label{7error}
\|x_k-y\|_{r,1}<\varkappa\,\varepsilon
\quad\mbox{при}\quad k\geq\varrho.
\end{equation}
Число $\varkappa$ можно выбрать независимым от $\varepsilon$, $\widehat{\varrho}$, $f$, $q$, $f_k$ и $q_k$, а число $\varrho$ можно выбрать независимым $f_k$ и $q_k$. В частности, если $f_k\to f$ в $(L_{1})^{m}$ и $q_k\to q$ в $\mathbb{C}^{rm}$ при $k\to\infty$, то $x_k\to y$ в $(W^r_1)^{m}$ при $k\to\infty$.\rm

При некотором более слабом условии, чем $\|f_k-f\|_{1}<\varepsilon$, можно вместо \eqref{7error} получить оценку погрешности $x_k-y$ в норме $\|\cdot\|_{(r-1)}$ (которая в свою очередь слабее соболевской нормы $\|\cdot\|_{r,1}$). Обозначим через $F$ і $F_k$ первообразные функций $f$ и $f_k$ на $[a,b]$, подчиненные условиям $F(a)=0$ и $F_k(a)=0$.

\textbf{Теорема 3.} \it Пусть заданы натуральное число $\widehat{\varrho}\geq\widetilde{\varrho}$ и действительное число $\varepsilon>0$. Предположим, что правые части краевых задач \eqref{Lk=fk}, \eqref{Bk=gk}
удовлетворяют условиям
\begin{equation}\label{cond-th3}
\|F_k-F\|_{C}<\varepsilon\quad\mbox{и}\quad\|q_k-q\|<\varepsilon
\quad\mbox{при}\;\;k\geq\widehat{\varrho}.
\end{equation}
Кроме того, предположим, что
\begin{equation}\label{supBk}
\sigma:=\sup
\bigl\{\|B_{k}:(C^{(r-1)})^{m}\to\mathbb{C}^{rm}\|:
k\geq\widehat{\varrho}\bigr\}<\infty.
\end{equation}
Тогда существуют положительные числа $\varkappa$ и $\varrho\geq\widehat{\varrho}$ такие, что
\begin{equation}\label{error(r-1)}
\|x_k-y\|_{(r-1)}<\varkappa\,\sigma\,\varepsilon
\quad\mbox{при}\quad k\geq\varrho.
\end{equation}
Число $\varkappa$ можно выбрать независимым от $\varepsilon$, $\sigma$, $\widehat{\varrho}$, $f$, $q$ и задач \eqref{Lk=fk}, \eqref{Bk=gk} (тогда $\varkappa$ будет зависеть только от $L$ и $B$), а число $\varrho$ можно выбрать независимым от $f_k$ и $q_k$. В~частности, если $F_k\to F$ в $(C)^{m}$ и $q_k\to q$ в $\mathbb{C}^{rm}$ при $k\to\infty$, то $x_k\to y$ в $(C^{(r-1)})^{m}$ при $k\to\infty$. \rm

Условие \eqref{cond-th3} слабее условия \eqref{74aut.r.s.est}, поскольку $\|F_{k}-F\|_{C}\leq\|f_{k}-f\|_{1}$.
Заметим, что условию \eqref{supBk} удовлетворяет построенная в доказательстве теоремы~1 аппроксимирующая последовательность многоточечных краевых задач вида \eqref{73MSp.syste}, \eqref{7MSp.2B_m}.

\textbf{3. Обоснование результатов.} Обозначим через $\mathrm{NBV}$ комплексное линейное пространство всех функций $\nobreak{g:[a,b]\to\mathbb{C}}$ ограниченной вариации на $[a,b]$, которые непрерывны слева на $(a,b)$ и удовлетворяют условию $g(a)=0$. Это пространство снабжено нормой, равной полной вариации $\operatorname{V}(g,[a,b])$ функции $g$ на $[a,b]$. Оно является полным, несепарабельным и нерефлексивным относительно этой нормы \cite[разд.~IV, п. 12 и~15]{DanfordShvarts1958}.

По теореме Ф.~Риса, каждый непрерывный линейный функционал $\ell$ на комплексном банаховом пространстве $C$ представляется единственным образом в виде
\begin{equation}\label{7l}
\ell(x)=\int\limits_a^b x(t)d\varphi(t)\quad\mbox{для любого}\;\;x\in C,
\end{equation}
где $\varphi$~--- некоторая функция класса $\mathrm{NBV}$. Более того, соответствующее отображение $\mathcal{I}:\nobreak \varphi \mapsto \ell$ порождает изометрический изоморфизм $\mathcal{I} \colon \mathrm{NBV} \leftrightarrow C^\prime$ (см., например, \cite[гл.~IV, п.~13, упражнение~35] {DanfordShvarts1958}). Здесь, как обычно, $C^\prime$~--- банахово пространство, сопряженное к $C$ и наделенное нормой функционалов. Интеграл в \eqref{7l} понимается как интеграл Римана\,--\,Стилтьеса или Лебега\,--\,Стилтьеса по комплексной мере, порожденной функцией~$\varphi$.

Пространство $\mathrm{NBV}$ изометрически изоморфно банахову пространству всех комплексных борелевых мер на $[a, b]$ (норма в последнем равна полной вариации меры) \cite[гл.~IV, п.~12]{DanfordShvarts1958}. Простейшим примером таких мер служат меры Дирака с одноточечными носителями. Обозначим через $\chi_{\gamma}$ характеристическую функцию (индикатор) подмножества $\gamma$ отрезка $[a, b]$. При этом изоморфизме таким мерам Дирака соответствуют характеристические функции $\chi_{(c, b]}$, где $a \leq c <b$, и $\chi_{\{b \}}$. Обозначим через $S$ комплексную линейную оболочку всех этих характеристических функций. В~доказательстве теоремы~1 ключевую роль будет играть такой результат:

\textbf{Утверждение 1.} \it
Множество $\mathcal{I} (S)$ является секвенциально плотным в пространстве $C^\prime $ в $w^*$-топологии, т.~е. для каждого функционала $\ell \in C^\prime$ существует последовательность $(\varphi_k)^\infty_{k=1} \subset S$ такая, что
\begin{equation} \label{74aut.t3}
   \lim_{k \to \infty} \int \limits_a^b x (t) d\varphi_k (t)=\ell (x) \quad \mbox{для любого}\;\;x\in C.
\end{equation} \rm

Как известно (см., например, \cite[п.~IV.5, с.~114]{ReedSimon}), множество $\mathcal{I} (S)$ плотно в пространстве $C^\prime$ в $w^*$-топологии. Однако, поскольку указанная топология неметризируема, то из этой плотности не следует утверждение~1. Его конструктивное доказательство приведено в \cite[теоремы 2 и 3]{MasliukPelekhataSoldatov2020MFAT}, где предложен явный алгоритм построения аппроксимирующей последовательности функций $\varphi_k \in S$.

\textbf{\emph{Доказательство теоремы} 1.} Поскольку по предположению множество $\mathcal{X}$ плотно в пространстве $(L_1)^{m \times m}$, выберем последовательности $(A_{r-l, k})_{k=1}^{\infty} \subset \mathcal{X}$, где $1 \leq l \leq r$, удовлетворяющие условиям
\begin{equation} \label{74aut.eq.1}
    \mathcal{X} \ni A_{r-l,k}\to A_{r-l} \quad\mbox{в}\quad (L_1)^{m \times m}\quad\mbox{при}\quad k \to \infty.
\end{equation}
Построим последовательность операторов $B_k$ вида \eqref{7MSp.2B_m} такую, что
\begin{equation} \label{74aut.eq.2}
    B_k\,y\to By \quad \mbox{в} \quad \mathbb{C}^{rm} \quad \mbox{для каждого} \quad y \in (C^{(r-1)})^{m}.
\end{equation}
Для этого используем однозначное представление непрерывного линейного оператора \eqref{72MSp.oper} в виде
\begin{equation} \label{7MSp.1predst}
By=\sum_{l=0}^{r-2}\alpha_{l}\,y^{(l)}(a)+
\int\limits_a^b(d\Phi(t))\,y^{(r-1)}(t)
\quad\mbox{для всех}\quad y\in(C^{(r-1)})^{m}.
\end{equation}
Здесь каждое $\alpha_l$~--- некоторая комплексная числовая матрица размера $rm \times m$, а
$$
\Phi (t)=(\varphi^{\lambda, \mu} (t))_{\substack {\lambda=1, \ldots, rm \\\mu=1, \ldots, m}}
$$
--- некоторая матрица-функция этого же размера такая, что каждое $\varphi^{\lambda,\mu}\in\mathrm{NBV}$; при этом интеграл понимается в смысле Римана\,--\,Стилтьеса. (Разумеется, если $r=1$, то в \eqref{7MSp.1predst} отсутствует сумма по индексу $l$.) Указанное представление следует непосредственно из известного описания пространства, сопряженного к $C^{(r-1)}$ (см., например, \cite[с.~344]{DanfordShvarts1958}).

Зафиксировав $\lambda\in\{1,\ldots,rm\}$ и $\mu\in\{1,\ldots,m\}$, определим функционал $\ell\in C^\prime$ по формуле \eqref{7l}, в которой положим $\varphi(t)\equiv\varphi^{\lambda,\mu}(t)$. Согласно утверждению~1 существует последовательность $(\varphi_{k}^{\lambda, \mu})^\infty_{k=1}\subset S$, удовлетворяющая условию \eqref{74aut.t3}, где $\varphi_k(t)\equiv\varphi_{k}^{\lambda,\mu}(t)$. Следовательно,
\begin{equation}\label{74uat.ell.lim}
    \lim_{k \to \infty}\int\limits_a^bz^{(r-1)}(t)d\varphi_k^{\lambda,\mu}(t)=\int \limits_a^bz^{( r-1)}(t)d\varphi^{\lambda,\mu}(t) \quad\mbox{для всех}\quad z\in C^{(r-1)}.
\end{equation}
Поскольку $\varphi_{k}^{\lambda, \mu} \in S$, то левый интеграл является линейной комбинацией значений функции $z^{(r-1)} (t)$ в некоторых точках отрезка $[a, b]$. Поэтому, полагая
\begin{equation*}
B_k\,y=\sum_{l=0}^{r-2}\alpha_{l}\,y^{(l)}(a)+
\int\limits_{a}^{b}\left(d(\varphi_k^{\lambda,\mu}(t))_{\substack{\lambda=1,
\ldots,rm \\\mu=1, \ldots, m}}\right) y^{(r-1)}(t)
\end{equation*}
для произвольного $y \in (C^{(r-1)})^{m}$, замечаем, что линейные непрерывные операторы $B_k:\nobreak (C^{(r-1)})^{m} \rightarrow \mathbb C^{rm}$, где $k \in \mathbb{N}$, принимают вид \eqref{7MSp.2B_m} и обладают нужным свойством \eqref{74aut.eq.2} ввиду \eqref{7MSp.1predst} и~\eqref{74uat.ell.lim}.

Обозначим через $\bar {B}$ и $\bar {B}_k$ соответственно сужение операторов $B$ и $B_k$ на пространство Соболева $(W^r_1)^m$. Поскольку последнее непрерывно вложено в пространство $(C^{(r-1)})^{m}$, то линейные операторы $\bar {B}$ и $\bar {B}_k$ непрерывно действуют из $(W^ r_1)^m$ в $\mathbb{C}^{rm}$. Рассмотрим задачу, состоящую из дифференциальной системы \eqref{72MSp.syste} и краевого условия
\begin{equation}\label{72MSp.ku.por.gen}
\bar{B}y=q.
\end{equation}
Кроме того, для каждого $k \in \mathbb{N}$ рассмотрим задачу, состоящую из дифференциальной системы \eqref{73MSp.syste} и краевого условия
\begin{equation}\label{72MSp.ku.por.mp}
\bar{B}_{k}\,y_{k}=q.
\end{equation}
Решения этих задач берутся в классе $(W^r_1)^m$. По предположению первая из них однозначно разрешима в пространстве $(W^r_1)^m$. Согласно \eqref{74aut.eq.2} имеем сходимость
\begin{equation} \label{74aut.eq.2.ner}
\bar{B}_k\,y\to\bar{B} y \quad \mbox{в} \quad \mathbb{C}^{rm} \quad
\mbox{для каждого}\quad y\in(W^{r}_1)^{m}.
\end{equation}
Предельные свойства этих задач исследованы в статьях \cite{GnypKodlyukMikhailets2015UMJ, GnypMikhailetsMurach2017}. На основании \cite[теорема~3]{GnypKodlyukMikhailets2015UMJ} (или \cite[теорема~2.3]{GnypMikhailetsMurach2017}) делаем вывод, что из условий \eqref{74aut.eq.1} и \eqref{74aut.eq.2.ner} следует однозначная разрешимость задачи \eqref{73MSp.syste}, \eqref{72MSp.ku.por.mp} (а следовательно, и задачи \eqref{73MSp.syste}, \eqref{7MSp.2B_m}) при $k \gg1$, и нужное асимптотическое свойство \eqref{7MSp.2granum}.

Остается отметить, что функции $\varphi_k^{\lambda, \mu} (t)$, а поэтому и операторы $B_k$ не зависят от $f$ и $q$, и строятся явно, как это было показано в \cite[доказательство лем~1 и 2]{MasliukPelekhataSoldatov2020MFAT}. Это касается и матриц $A_{r-l, k}$.

Теорема~1 доказана.

\textbf{\emph{Доказательство теоремы} 2.} Рассмотрим краевые задачи \eqref{72MSp.syste}, \eqref{72MSp.ku.por.gen} и \eqref{73MSp.syste}, \eqref{72MSp.ku.por.mp}, где (как и раньше) $\bar {B}$ и $\bar {B}_k$ обозначают соответственно сужение операторов $B$ и $B_k$ на пространство Соболева $(W^r_1)^m$. Свяжем с этими задачами линейные операторы $(L, \bar {B})$ и $(L_k, \bar {B}_k)$, которые непрерывно действуют из $(W^{r}_1)^{m}$ в $(L_1)^{m} \times \mathbb{C}^{rm}$. Эти операторы обратимы при $k \geq \widetilde{\varrho}$. Рассмотрим обратные к ним операторы $(L, \bar {B})^{-1}$ и $(L_k, \bar {B}_k)^{-1}$. Поскольку по предположению краевые задачи \eqref{73MSp.syste}, \eqref{7MSp.2B_m} удовлетворяют заключению теоремы~1, то
\begin{equation}\label{7f16}
(L_k,\bar{B}_k)^{-1}(f,q)=y_{k}\to y=(L,\bar{B})^{-1}(f,q)
\quad\mbox{в}\quad(W^{r}_1)^{m}\quad\mbox{при}\quad k \to \infty
\end{equation}
для произвольных $f \in (L_1)^{m}$ и $q \in \mathbb{C}^{rm}$, т.~е. $(L_k, \bar {B}_k)^{-1}$ сходится к $(L,\bar{B})^{-1}$ в сильной операторной топологии. Значит, по теореме Банаха\,--\,Штейнгауса
\begin{equation*}
\widetilde{\varkappa}=\sup \bigl \{\| (L_k, \bar {B}_k)^{-1}:
(L_1)^{m} \times \mathbb{C}^{rm} \to (W^{r}_1)^{m} \|:
k \geq \widetilde{\varrho} \bigr\} <\infty.
\end{equation*}
Следовательно,
\begin{gather*}
\|x_{k}-y\|_{r, 1} =
\|(L_k,\bar{B}_k)^{-1}(f_k,q_k)-(L,\bar{B})^{-1}(f,q)\|_{r,1}\leq\\
\leq \| (L_k, \bar {B}_k)^{- 1} (f_k, q_k) - (L_k, \bar {B}_k)^{-1}(f,q)\|_{r,1}+\\
\quad+\|(L_k,\bar{B}_k)^{-1}(f, q)-(L,\bar{B})^{-1}(f,q)\|_{r,1}\leq\\
\leq \widetilde{\varkappa} (\| f_k-f \|_{1} + \|q_k-q\|)+\| y_{k}-y\|_{r,1}<\\
<2\widetilde{\varkappa}\,\varepsilon+\varepsilon=\varkappa\,\varepsilon
\quad \mbox{при} \quad k \geq \varrho \geq \widehat {\varrho}
\end{gather*}
на основании \eqref{74aut.r.s.est} и \eqref{7f16}. Здесь $\varkappa:=2\widetilde{\varkappa}+1$, а число $\varrho\geq\widehat{\varrho}$ удовлетворяет импликации $k \geq \varrho \Rightarrow \| y_{k}-y\|_{r,1}<\varepsilon$. При этом $\varrho$ не зависит от правых частей $f_k$ и $q_k$, поскольку от них не зависит $y_{k} -y$. Таким образом, решение $x_k$ задачи \eqref{Lk=fk}, \eqref{Bk=gk} обладает нужным свойством~\eqref{7error}. Теорема~2 доказана.

\textbf{\emph{Доказательство теоремы} 3.} Все границы рассматриваем при $k \to \infty$. Исследуем сначала случай $r=1$. Если $f_k \to f$ в $(L_{1})^{m}$ и $q_k \to q$ в $\mathbb{C}^{m}$, то $x_k \to y$ в $( W^1_1)^{m}$ согласно теореме~2. Следовательно, рассматривая $y$ и $x_{k}$ как решения краевых задач \eqref{72MSp.syste}, \eqref{72MSp.ku.por.gen} и \eqref{Lk=fk}, $\bar{B}_{k}\,x_{k}=q_{k}$, заключаем на основании \cite[теорема~2.3]{GnypMikhailetsMurach2017}, что
\begin{gather} \label{Ak-to-A0}
A_{0,k}\to A_{0}\;\;\mbox{в}\;\;(L_{1})^{m \times m},\\
B_{k}x\to Bx\;\;\mbox{в}\;\;\mathbb{C}^{m}\;\;
\mbox{для каждого}\;\;x\in(W^{1}_{1})^{m}.\label{Bk-to-B}
\end{gather}
Пусть натуральное $k \geq \widehat {\varrho}$. Напомним, что $y_{k}$~--- единственное решение задачи \eqref{73MSp.syste}, \eqref{7MSp.2B_m}. Имеем:
\begin{equation} \label{triangle-inequal}
\|x_k-y\|_{C}\leq\|x_k-y_k\|_{C}+\|y_k-y\|_{C},
\quad\mbox{где}\quad\|y_k-y\|_{C}\to0
\end{equation}
согласно теореме~1 и непрерывному вложению $(W^{1}_{1})^{m}\hookrightarrow(C)^{m}$. Оценим норму $\|x_k-y_k\|_{C}$.

Обозначим через $Y_{k}$ матрицант системы \eqref{73MSp.syste}, отнесенный к точке $t=a$, т.~е. матрица-функция $Y_{k} \in (W^{1}_ {1})^{m \times m}$ является решением задачи Коши
\begin{equation}
Y_{k}'(t)=-A_{0,k}(t)Y_{k}(t)\;\;\mbox{для п.в.}\;\;t\in[a,b],\quad\quad Y_{k}(a)=E_{m},
\end{equation}
где $E_{m}$~--- единичная матрица порядка~$m$. Как хорошо известно, числовая матрица $Y_{k}(t)$ невырождена для каждого $t \in [a, b]$. Обозначим через $[B_{k} Y_{k}]$ квадратную числовую матрицу порядка $m$, произвольный $j$-й столбец которой является результатом действия оператора $B_{k}$ на $j$-й столбец матрицанта~$Y_{k}$. Поскольку задача \eqref{Lk=fk}, \eqref{Bk=gk} однозначно разрешима, то матрица $[B_{k} Y_{k}]$ невырождена согласно \cite[формула (1.7)]{Kiguradze1987}.

Как известно \cite[формула (1.26)]{Kiguradze1987}, решение $x_k$ краевой задачи \eqref{Lk=fk}, \eqref{Bk=gk} представляется в виде
\begin{equation} \label{hat-yk-represent}
x_{k}(t)=z_{k}(t,q_k)+Y_{k}(t)h_{k}(f_k)+R_{k}(t,f_k)
\;\;\mbox{для всех}\;\;t\in[a, b].
\end{equation}
Здесь
\begin{gather} \label{zk-def}
z_{k}(t,q_k):=Y_{k}(t)[B_{k}Y_{k}]^{-1}q_k,\\
h_{k}(f_k):=-[B_{k}Y_{k}]^{-1}B_{k}R_{k}(\cdot,f_k), \label{hk-def}
\end{gather}
где
\begin{equation} \label{Rk-def}
R_{k}(t,f_k):=Y_{k}(t)\int\limits_{a}^{t} Y_{k}^{-1}(\tau)f_k(\tau)d\tau=
F_{k}(t)-Y_{k}(t)
\int\limits_{a}^{t}Y_{k}^{-1}(\tau)A_{0, k}(\tau)F_k(\tau)d\tau.
\end{equation}
Взяв $f_k=f$ и $q_k=q$ в формулах \eqref{zk-def}--\eqref{Rk-def}, получим разложение
\begin{equation} \label{yk-represent}
y_{k}(t)=z_{k}(t,q)+Y_{k}(t)h_{k}(f)+R_{k}(t,f)
\;\;\mbox{для всех}\;\;t\in[a,b].
\end{equation}

Обозначим через $Y$ матрицант системы \eqref{72MSp.syste} отнесенный к точке $t=a$. Из свойства \eqref{Ak-to-A0} следует в силу \cite[теорема~2.3]{GnypMikhailetsMurach2017}, что $Y_{k} \to Y$ в $(W^{1}_{1})^{m \times m} \hookrightarrow (C)^{m \times m}$
(вложение непрерывно). Поэтому $Y_{k}^{- 1} \to Y^{- 1}$ в $(C)^{m \times m}$.
Кроме того,
\begin{equation*}
[B_{k}Y_{k}]=[B_{k}(Y_{k}-Y)]+[B_{k}Y]\to[BY]\;\;
\mbox{в}\;\;\mathbb{C}^{m \times m}
\end{equation*}
согласно \eqref{supBk} и \eqref{Bk-to-B}. Матрица $[BY]$ невырождена, поскольку задача \eqref{72MSp.syste}, \eqref{72MSp.kue} однозначно разрешима. Поэтому $[B_{k} Y_{k}]^{- 1} \to [BY]^{- 1}$ в $\mathbb{C}^{m \times m}$. Следовательно, существует целое число $\varrho_{1} \geq \widehat {\varrho}$ такое, что
\begin{equation}\label{bounds}
\|Y_{k}\|_{C}\cdot\|[B_{k}Y_{k}]^{-1}\|\leq 1+\|Y\|_{C}\cdot\|[BY]^{-1}\|=:c_{1}\;\;\mbox{при}\;\,k\geq\varrho_{1}
\end{equation}
и, кроме того,
\begin{equation}\label{Ak-bound}
\|Y_{k}\|_{C}\cdot\|Y_{k}^{-1}\|_{C}\cdot\|A_{0,k}\|_{1}\leq
1+\|Y\|_{C}\cdot\|Y^{-1}\|_{C}\cdot\|A_{0}\|_{1}=:c_{2}-1
\;\;\mbox{при}\;\;k\geq\varrho_{1}
\end{equation}
в силу \eqref{Ak-to-A0}. Число $\varrho_{1}$ не зависит от правых частей задач \eqref{72MSp.syste}, \eqref{72MSp.kue} и \eqref{Lk=fk}, \eqref{Bk=gk}.

Воспользовавшись условиями \eqref{cond-th3}, \eqref{supBk} и свойствми \eqref{bounds}, \eqref{Ak-bound}, оценим разности соответствующих слагаемых в разложениях \eqref{hat-yk-represent} и \eqref{yk-represent}. Пусть натуральное $k \geq \varrho_{1}$. Тогда
\begin{equation*}
\|z_{k}(\cdot,q_k)-z_{k}(\cdot,q)\|_{C}=
\|Y_{k}(\cdot)[B_{k}Y_{k}]^{-1}(q_k-q)\|_{C}\leq
c_{1}\|q_k-q\|<c_{1}\varepsilon.
\end{equation*}
Кроме того,
\begin{equation}\label{Rk-bound}
\begin{gathered}
\|R_{k}(\cdot,f_k)-R_{k}(\cdot,f)\|_{C}
\leq\|F_{k}-F\|_{C}+\|Y_{k}\|_{C}\,
\|Y_{k}^{-1}A_{0,k}(F_k-F)\|_{1}\leq\\
\leq\|F_{k}-F\|_{C}+
\|Y_{k}\|_{C}\,\|Y_{k}^{-1}\|_{C}\,\|A_{0,k}\|_{1}\,\|F_k-F\|_{C}<
c_{2}\varepsilon.
\end{gathered}
\end{equation}
Поэтому
\begin{equation} \label{Ykhk-bound}
\begin{gathered}
\|Y_{k}(\cdot)h_{k}(f_k)-Y_{k}(\cdot)h_{k}(f)\|_{C}=
\|Y_{k}(\cdot)[B_{k}Y_{k}]^{-1}B_{k}
(R_{k}(\cdot,f_k)-R_{k}(\cdot,f))\|_{C}\leq\\
\leq\|Y_{k}\|_{C}\,\|[B_{k}Y_{k}]^{-1}\|\cdot\sigma\,
\|R_{k}(\cdot,f_k)-R_{k}(\cdot,f)\|_{C}<
c_{1}c_{2}\,\sigma\,\varepsilon.
\end{gathered}
\end{equation}

Таким образом,
\begin{gather*}
\|x_{k}-y_{k}\|_{C}\leq (c_{1}+c_{2}+c_{1}c_{2}\,\sigma)\varepsilon=
((c_{1}+c_{2})\sigma^{-1}+c_{1}c_{2})\,\sigma\,\varepsilon\leq
\varkappa_{0}\sigma\,\varepsilon\quad\mbox{при}\;\;k\geq\varrho_{1},
\end{gather*}
где число
\begin{equation*}
\varkappa_{0}:=
(c_{1}+c_{2})\,\|B:(C)^{m}\to\mathbb{C}^{m}\|^{-1}+c_{1}c_{2}
\end{equation*}
зависит только от $A_{0}$ и $B$. Отсюда на основании \eqref{triangle-inequal} заключаем, что $\|x_{k}-y\|_{C}<(1+\varkappa_{0}) \sigma\varepsilon$ при $k\geq\varrho$, где число $\varrho\geq\varrho_1$ удовлетворяет импликации $k \geq\varrho\Rightarrow \|y_{k}-y\|_{C}<\sigma\varepsilon$ и не зависит от $f_k$ и $q_k$. Остается положить $\varkappa:=1+\varkappa_{0}$. Теорема~3 доказана в случае $r=1$.

Перейдем к случаю $r \geq2$. Сведем краевые задачи \eqref{72MSp.syste}, \eqref{72MSp.kue} и \eqref{Lk=fk}, \eqref{Bk=gk} к краевым задачам для систем дифференциальных уравнений первого порядка. Начнем с первой из них. Как обычно, положим $g:=\mathrm{col}(0,f)\in(L_1)^{rm}$ и
\begin{equation} \label{P-def}
P:=\left(
\begin{array}{ccccc}
O_m & -E_m & O_m & \ldots & O_m \\
O_m & O_m & -E_m & \ldots & O_m \\
\vdots & \vdots & \vdots & \ddots & \vdots \\
O_m & O_m & O_m & \ldots & -E_m \\
A_0 & A_1 & A_2 & \ldots & A_{r-1}\\
\end{array}\right)\in(L_1)^{rm\times rm},
\end{equation}
где $O_m$~--- квадратная нуль-матрица порядка $m$. Воспользовавшись представлением краевого оператора $B$ в виде \eqref{7MSp.1predst}, положим
\begin{equation} \label{T-def}
Tv:=
\sum_{l=0}^{r-2}\alpha_{l}\,v^{l}(a)+\int\limits_a^b(d\Phi(t))\,v^{r-1}(t)
\quad \mbox{для произвольного} \quad v \in (C)^{rm},
\end{equation}
где $v=\mathrm{col} (v^{0}, \ldots, v^{r-1})$ и $v^{l} \in (C)^{m}$ для каждого $l \in \{0, \ldots, r-1 \}$. Имеем ограниченный линейный оператор $T: (C)^{rm} \to \mathbb{C}^{rm}$. Если $y \in (W_1^{r})^m$ является решением задачи \eqref{72MSp.syste}, \eqref{72MSp.kue}, то вектор-функция
\begin{equation} \label{reduced-solution}
v=\mathrm{col}\bigl(y,y',\ldots,y^{(r-1)}\bigr)
\in(W^{1}_{1})^{rm}
\end{equation}
является решением задачи
\begin{align}\label{reduced-system}
v'(t)+P(t)v(t)&=g(t)
\quad\mbox{для п.в.}\;\;t\in[a,b],\\
Tv&=q. \label{reduced-bound-cond}
\end{align}
И наоборот, если $v \in (W^{1}_{1})^{rm}$ является решением задачи \eqref{reduced-system}, \eqref{reduced-bound-cond}, то вектор-функция
$y:=v^{0}$ принадлежит пространству $(W_1^{r})^m$ и является решением задачи \eqref{72MSp.syste}, \eqref{72MSp.kue}, причем выполняется равенство~\eqref{reduced-solution}.

Аналогично сведем каждую задачу \eqref{Lk=fk}, \eqref{Bk=gk} к краевой задаче
\begin{align}\label{reduced-system-k}
u_{k}'(t)+P_{k}(t)u_{k}(t)&=
g_{k}(t)\quad\mbox{для п.в.}\;\;t\in[a,b],\\
T_{k}u_{k}&=q_{k}. \label{reduced-bound-cond-k}
\end{align}
Здесь решение $u_{k} \in (W_1^{1})^{rm}$, $g_{k}:=\mathrm{col} (0, f_{k}) \in (L_1)^{rm}$, матрица-функция $P_{k} \in (L_1)^{rm \times rm}$ определена по формуле \eqref{P-def}, в которой каждое $A_{j}$ заменено на $A_{j, k}$. Ограниченный линейный оператор $T_{k}:(C)^{rm}\to\mathbb{C}^{rm}$ определен подобно оператору $T$ с помощью однозначного представления оператора $B_{k}$ в виде
\begin{equation} \label{Bk-reprent}
B_{k}x=\sum_{l=0}^{r-2}\alpha_{l,k}\,x^{(l)}(a)+
\int\limits_a^b(d\Phi_{k}(t))\,x^{(r-1)}(t)
\quad \mbox{для произвольного} \quad x \in (C^{(r-1)})^{m},
\end{equation}
где каждое $\alpha_{l, k}$~--- некоторая комплексная числовая матрица размера $rm \times m$, а $\Phi_{k} (t)$ --- некоторая матрица-функция этого же размера такая, что каждый ее элемент принадлежит $\mathrm{NBV}$. А именно:
\begin{equation} \label{Tk-def}
T_{k}u:=
\sum_{l=0}^{r-2}\alpha_{l,k}\,u^{l}(a)+
\int\limits_a^b(d\Phi_{k}(t))\,u^{r-1}(t)
\quad \mbox{для произвольного} \quad u \in (C)^{rm},
\end{equation}
где $u=\mathrm{col} (u^{0}, \ldots, u^{r-1})$ и $u^{l} \in (C)^{m}$ для каждого $l \in \{0, \ldots, r-1 \}$.

Если $f_k \to f$ в $(L_{1})^{m}$ и $q_k \to q$ в $\mathbb{C}^{rm}$, то $x_k \to y$ в $( W^r_1)^{m}$ согласно теореме~2. Поэтому на основании \cite[теорема~2.3]{GnypMikhailetsMurach2017} имеем: $A_{r-l,k} \to A_{r-l}$ в $(L_{1})^{m \times m}$ для каждого $l\in\{1,\ldots,r\}$ и $B_{k}x\to Bx$ в $\mathbb{C}^{rm}$ для каждого $x\in(W^{r}_{1})^{m}$. Итак,
\begin{equation} \label{Pk-to-P}
P_{k}\to P\;\;\mbox{в}\;\;(L_{1})^{rm \times rm}
\end{equation}
и, кроме того,
\begin{equation} \label{Bk-to-B-strong}
B_{k}x\to Bx\;\;\mbox{в}\;\;\mathbb{C}^{rm}\;\;
\mbox{для каждого}\;\;x\in(C^{(r-1)})^{m},
\end{equation}
учитывая условие \eqref{supBk} и плотность множества $(W^{r}_{1})^{m}$ в пространстве $(C^{(r-1)})^{m}$. Ввиду представлений \eqref{7MSp.1predst} и \eqref{Bk-reprent} свойство \eqref{Bk-to-B-strong} эквивалентно системе следующих четырех условий:
\begin{itemize}
  \item[(i)] $\alpha_{l,k} \to \alpha_{l}$ в $\mathbb{C}^{rm \times m}$ для каждого $l \in \{0,\ldots,r-2\}$;
  \item[(ii)] $\sup \{\|V (\Phi_k [a, b]) \|: k \geq \widetilde{\varrho}\,\}
      <\infty$;
  \item[(iii)] $\Phi_k (b) \to \Phi (b)$ в $\mathbb{C}^{rm \times m}$;
  \item[(iv)] $\int_{a}^{t} \Phi_k (s)ds\to
\int_{a}^{t} \Phi (s) ds$ в $\mathbb{C}^{rm \times m}$ для каждого $t \in (a, b]$.
\end{itemize}
Здесь числовая матрица $V(\Phi_k, [a,b])$ образована полными вариациями элементов матрицы-функции~$\Phi_k$. Эта эквивалентность следует из критерия Ф.~Рисса слабой сходимости непрерывных линейных функционалов на пространстве $C$ (см., например, \cite[гл.~III, п.~55]{RieszSz-Nagy90}). Следовательно, \eqref{Bk-to-B-strong} влечет сходимость
\begin{equation} \label{Tk-to-T-strong}
T_{k}u\to Tu\;\;\mbox{в}\;\;\mathbb{C}^{rm}\;\;
\mbox{для каждого}\;\;u\in(C)^{rm}.
\end{equation}

В этом абзаце рассмотрим краевые задачи \eqref{reduced-system}, \eqref{reduced-bound-cond} и \eqref{reduced-system-k}, \eqref{reduced-bound-cond-k}, где $k \geq \widetilde{\varrho}$, для \emph{произвольных} правых частей $g, g_k \in (L_1)^{rm}$ и $q, q_k \in \mathbb{C}^{rm}$. Они принадлежат к классу задач, исследованных в этой работе в случае $r=1$. Они однозначно разрешимы в пространстве $(W^{1}_{1})^{rm}$, поскольку имеют лишь тривиальное решение в случае нулевых правых частей, что следует из наших предположений об однозначной разрешимости краевых задач \eqref{72MSp.syste}, \eqref{72MSp.kue} и
\eqref{Lk=fk}, \eqref{Bk=gk} в пространстве $(W^{r}_{1})^{m}$. Из свойств \eqref{Pk-to-P} и \eqref{Tk-to-T-strong} следует согласно \cite[теорема~2.3]{GnypMikhailetsMurach2017}, что
\begin{equation*}
\bigl(g_k\to g\;\,\mbox{в}\;\,(L_1)^{rm}\bigr)\,\wedge\,
\bigl(q_k\to q\;\,\mbox{в}\;\,\mathbb{C}^{rm}\bigr)\;\Longrightarrow\;
\bigl(u_k\to v\;\,\mbox{в}\;\,(W^{1}_{1})^{rm}\bigr).
\end{equation*}
В частности, последовательность задач \eqref{reduced-system-k}, \eqref{reduced-bound-cond-k} удовлетворяет заключению теоремы~1, рассмотренной для задачи \eqref{reduced-system}, \eqref{reduced-bound-cond}.

Из формул \eqref{Bk-reprent} і \eqref{Tk-def} вытекает, что
\begin{equation*}
\|B_{k}:(C^{(r-1)})^{m}\to\mathbb{C}^{rm}\|\leq
\|T_{k}:(C)^{rm}\to\mathbb{C}^{rm}\|\leq
\theta\,\|B_{k}:(C^{(r-1)})^{m}\to\mathbb{C}^{rm}\|,
\end{equation*}
где $\theta$~--- некоторое положительное число, которое может зависеть лишь от $r$, $m$ и $b-a$. Следовательно, условие \eqref{supBk} влечет за собой неравенство
\begin{equation}\label{norm-Tk-bound}
\|T_{k}:(C)^{rm}\to\mathbb{C}^{rm}\|\leq \theta\sigma
\quad\mbox{при}\;\;k\geq\widehat{\varrho}.
\end{equation}

Пусть снова $g=\mathrm{col} (0, f)$ и $g_{k}=\mathrm{col} (0, f_{k})$. Тогда
\begin{equation} \label{v-uk-special}
v=\mathrm{col}\bigl(y,y',\ldots,y^{(r-1)})\quad\mbox{и}\quad
u_k=\mathrm{col}\bigl(x_k,x_k',\ldots,x_k^{(r-1)})
\end{equation}
~--- решения задач \eqref{reduced-system}, \eqref{reduced-bound-cond} и \eqref{reduced-system-k}, \eqref{reduced-bound-cond-k} соответственно.
На основании теоремы~3, уже доказанной для этих задач, заключаем, что из условий \eqref{cond-th3} и \eqref{norm-Tk-bound} вытекает существование положительных чисел $\varkappa_{1}$ и $\varrho \geq \widehat {\varrho}$ таких, что
\begin{equation}\label{error-uk-v}
\|u_k-v\|_{C}<\varkappa_{1}\theta\sigma\varepsilon
\quad\mbox{при}\quad k\geq\varrho.
\end{equation}
Здесь $\varkappa_{1}$ зависит лишь от $P$ и $T$ (т.~е. только от $L$ и $B$), а число $\varrho$ можно выбрать независимым от $g_{k}$ (т.~е. от $f_k$) и $q_k$. Из формул \eqref{v-uk-special} и \eqref{error-uk-v} немедленно следует нужное неравенство \eqref{error(r-1)}, где $\varkappa:=\varkappa_{1}\theta$.

Теорема~3 доказана.

\textbf{4. Заключительные замечания.} Первые два из них касаются основного результата~--- теоремы~1.

\textbf{\textit {Замечание} 1.} В доказательстве теоремы~1 построена последовательность многоточечных краевых операторов $B_{k}$ вида \eqref{7MSp.2B_m}, которая сходится к оператору $B$ в топологии сильной сходимости непрерывных линейных операторов на паре пространств $(C^{(r-1)})^{m}$ и $\mathbb{C}^{rm}$. Для равномерной сходимости (т.~е. сходимости по операторной норме) такую аппроксимирующую последовательность построить, вообще говоря, нельзя. Действительно, ввиду представления \eqref{7MSp.1predst} оператора $B$, его аппроксимация операторами $B_{k}$ является по сути аппроксимацией подынтегральной функции $\Phi$ кусочно-постоянными функциями. Если, например, $r=1$, $m=1$ и $\Phi (t) \equiv t$, то равномерная сходимость $B_{k} \rightrightarrows B$ равносильна тому, что полная вариация $\operatorname {V} (\Phi_{k} - \Phi,[a, b]) \to0$ для некоторой последовательности кусочно-постоянных функций $\Phi_{k}$ на $[a, b]$. Для рассматриваемой функции $\Phi (t) \equiv t$ последняя сходимость невозможна, поскольку $\operatorname {V} (\Phi_{k}-\Phi, [a, b]) \geq b-a$.

\textbf{\textit {Замечание} 2.} Пусть целое $r \geq2$. Младшие части многоточечных краевых дифференциальных операторов $B_{k}$, построенных в доказательстве теоремы~1, не зависят от $k$ и принимают вид
\begin{equation*}
\sum_{l=0}^{r-2}\alpha_{l}\,y^{(l)}(a),
\end{equation*}
где числовые матрицы $\alpha_{l}$ взяты из представления краевого оператора $B$ в виде~\eqref{7MSp.1predst}.

\textbf{\textit {Замечание} 3.} Укажем допустимое значение постоянной $\varkappa$ в теореме~3. Как видно из ее доказательства, в случае $r=1$ можно взять
\begin{equation} \label{kappa-value}
\varkappa:=(c_{1}+c_{2})\lambda+c_{1}c_{2}+1,
\end{equation}
где
\begin{align*}
c_{1}&:=1+\|Y\|_{C}\cdot\|[BY]^{-1}\|,\\
c_{2}&:=2+\|Y\|_{C}\cdot\|Y^{-1}\|_{C}\cdot\|A_{0}\|_{1},\\
\lambda &:=\|B:(C)^{m}\to\mathbb{C}^{m}\|^{-1}.
\end{align*}
В случае $r \geq2$ постоянную $\varkappa$ можно определить по формуле \eqref{kappa-value}, где
\begin{align*}
c_{1}&:=1+\|V\|_{C}\cdot\|[BV^{\circ}]^{-1}\|,\\
c_{2}&:=2+\|V\|_{C}\cdot\|V^{-1}\|_{C}\,(b-a+\|A_{r-1}\|_{1}),\\
\lambda &:=\|B:(C^{(r-1)})^{m}\to\mathbb{C}^{rm}\|^{-1}.
\end{align*}
Здесь $V$~--- матрицант системы \eqref{reduced-system}, отнесенный к точке $t=a$ (т.~е. $V\in(W_{1}^{1})^{rm\times rm}$~--- решение задачи Коши $V'(t)=-P(t)V(t)$ для п.в. $t\in[a,b]$, $V(a)=E_{rm}$), $V^{\circ}$~--- матрица-функция размера $m \times rm$, образованная первыми $m$ строками матрицы-функции $V$ (все элементы матрицы $V^{\circ}$ принадлежат пространству $W_{1}^{r}$ ввиду \eqref{P-def}). Непосредственное использование доказательства теоремы~3 в случае $r \geq2$ приводит к более грубой оценке $\varkappa$ сверху. Однако, анализируя доказательство в случае $r=1$ применительно к приведенным краевым задачам \eqref{reduced-system}, \eqref{reduced-bound-cond} и \eqref{reduced-system-k}, \eqref{reduced-bound-cond-k}, приходим к указанному результату. При этом используем следующие обстоятельства:

1) Поскольку все столбцы матрицанта $V$ имеют вид \eqref{v-uk-special}, то $[BV^{\circ}]=[TV]$ для оператора $T$, определенного по формуле \eqref{T-def}, что и дает оценку $\| z_{k} (\cdot, q_k) -z_{k} (\cdot, q) \|_{C} \leq c_{1} \varepsilon$ при $k \gg1$.

2) Оценка \eqref{Rk-bound} для приведенных краевых задач принимает вид
\begin{gather*}
\|R_{k}(\cdot,g_k)-R_{k}(\cdot,g)\|_{C}
\leq\|G_{k}-G\|_{C}+\|V_{k}\|_{C}\,
\|V_{k}^{-1}P_{k}(G_k-G)\|_{1}\leq\\
\leq\|F_{k}-F\|_{C}+
\|V_{k}\|_{C}\,\|V_{k}^{-1}\|_{C}\,\|P_{k}(G_k-G)\|_{1}=\\
=\|F_{k}-F\|_{C}+
\|V_{k}\|_{C}\,\|V_{k}^{-1}\|_{C}\bigl(\|F-F_{k}\|_{1}+
\|A_{r-1,k}(F_{k}-F)\|_{1}\bigr)\leq\\
\leq\|F_{k}-F\|_{C}+
\|V_{k}\|_{C}\,\|V_{k}^{-1}\|_{C}\,(b-a+\|A_{r-1,k}\|_{1})\|F_{k}-F\|_{C}<
c_{2}\varepsilon
\end{gather*}
при $k\gg1$. Здесь $V_{k}$~--- матрицант системы \eqref{reduced-system-k}, отнесенный к точке $t=a$, и $G_k:=\mathrm{col}(0,F_k)$,  $G:=\mathrm{col}(0,F)$~--- вектор-функции размерности~$rm$.

3) Оценка \eqref{Ykhk-bound} для приведенных краевых задач принимает вид
\begin{gather*}
\|V_{k}(\cdot)h_{k}(g_k)-V_{k}(\cdot)h_{k}(g)\|_{C}=
\|V_{k}(\cdot)[T_{k}V_{k}]^{-1}T_{k}
(R_{k}(\cdot,g_k)-R_{k}(\cdot,g))\|_{C}=\\
=\|V_{k}(\cdot)[B_{k}V_{k}^{\circ}]^{-1}B_{k}
(R_{k}^{\circ}(\cdot,g_k)-R_{k}^{\circ}(\cdot,g))\|_{C}\leq\\
\leq\|V_{k}\|_{C}\,\|[B_{k}V_{k}^{\circ}]^{-1}\|\cdot\sigma\,
\|R_{k}^{\circ}(\cdot,g_k)-R_{k}^{\circ}(\cdot,g)\|_{(r-1)}<
c_{1}c_{2}\,\sigma\,\varepsilon
\end{gather*}
при $k\gg1$. Здесь матрица-функция $V_{k}^{\circ}$ образована первыми $m$ строками матрицы-функции $V_{k}$, вектор-функция $R_{k}^{\circ}(\cdot,\cdot)$ образована первыми $m$ компонентами вектор-функции $R_{k}(\cdot,\cdot)$, причем вектор-функции $R_{k}(\cdot,g_k)$ і $R_{k}(\cdot,g)$ принимают вид~\eqref{v-uk-special} в силу формул \eqref{hat-yk-represent} и \eqref{yk-represent} для приведенных задач.


\begin{thebibliography}{99}


\bibitem{Kiguradze1975}
\emph{Кигурадзе И. Т.} Некоторые сингулярные краевые задачи для обыкновенных дифференциальных уравнений.~-- Тбилиси: Изд-во Тбилисского ун-та, 1975.~-- 352~с.

\bibitem{Kiguradze1987}
\emph{Kiguradze I. T.} Boundary-value problems for systems of ordinary differential equations~// J.~Soviet Math.~-- 1988.~-- \textbf{43}.~-- P.  2259\,--\,2339.

\bibitem{Kiguradze2003}
\emph{Kiguradze I. T.} On boundary value problems for linear differential systems with singularities~// Differ. Equ.~-- 2003.~-- \textbf{39}, №~2.~-- P. 212\,--\,225.

\bibitem {Ashordia1996}
\emph{Ashordia M.} Criteria  of  correctness  of  linear  boundary  value  problems  for  systems of  generalized  ordinary  differential  equations~// Czechoslovak Math.~J.~-- 1996.~-- \textbf{46}, №~3.~-- P. 385\,--\,404.

\bibitem{MikhailetsReva2008DAN9}
\emph{Михайлец В.~А., Рева Н.~В.} Обобщения теоремы Кигурадзе о корректности линейных краевых задач~// Доп. НАН Укрїни.~-- 2008.~--  №~9.~-- С. 23\,--\,27.

\bibitem {KodliukMikhailetsReva2013}
\emph{Kodlyuk (Kodliuk) T.  I., Mikhailets V.  A., Reva N.  V.}  Limit  theorems  for  one-dimensional boundary-value problems~// Ukrainian Math.~J.~-- 2013.~-- \textbf{65}, №~1.~-- P. 77\,--\,90.

\bibitem{MikhailetsChekhanova2015JMathSci}
\emph{Mikhailets V.~A., Chekhanova G.~A.} Limit theorems for general one-dimensional boundary-value problems~// J.~Math. Sci. (N. Y.)~-- 2015.~-- \textbf{204}, №~3.~-- P. 333\,--\,342.

\bibitem{MikhailetsPelekhataReva2018UMJ}
\emph{Mikhailets V. A., Pelekhata O. B., Reva N. V.} Limit theorems for the solutions of boundary-value problems~// Ukrainian Math.~J.~-- 2018.~-- \textbf{70}, №~2.~-- P. 243\,--\,251.

\bibitem{PelekhataReva2019UMJ}
\emph{Pelekhata O. B., Reva N. V.} Limit theorems for the solutions of linear boundary-value problems for systems of differential equations~//
Ukrainian Math.~J.~-- 2019.~-- \textbf{71}, №~7.~-- P. 1061\,--\,1070.

\bibitem{MikhailetsReva2008DAN8}
\emph{Михайлец В. А., Рева Н.~В.} Предельный переход в системах линейных дифференциальных уравнений~// Доп. НАН Укрїни.~-- 2008.~-- №~8.~--  С. 28\,--\,30.



\bibitem{KodliukMikhailets2013JMS}
\emph{Kodliuk T. I., Mikhailets V. A.} Solutions of one-dimensional boundary-value problems with a parameter in Sobolev spaces~// J.~Math. Sci. (N.~Y.)~-- 2013.~-- \textbf{190}, №~4.~-- P. 589\,--\,599.


\bibitem{GnypKodlyukMikhailets2015UMJ}
\emph{Gnyp E. V., Kodlyuk (Kodliuk) T. I., Mikhailets V. A.} Fredholm boundary-value problems with parameter in Sobolev spaces~//  Ukrainian Math.~J.~-- 2015.~-- \textbf{67}, №~5.~-- P. 658\,--\,667.

\bibitem{Soldatov2015UMJ}
\emph{Soldatov V.~O.} On the continuity in a parameter for the solutions of boundary-value problems total with respect to the spaces  $C^{(n+r)}[a,b]$~// Ukrainian Math.~J.~--  2015.~-- \textbf{67}, №~5.~--  P. 785\,--\,794.

\bibitem{MikhailetsMurachSoldatov2016EJQTDE}
\emph{Mikhailets V. A., Murach A. A., Soldatov V.} Continuity in a parameter of solutions to generic boundary-value problems~// Electron.~J. Qual. Theory Differ. Equ.~--  2016.~-- №~87.~-- P. 1\,--\,16.

\bibitem{MikhailetsMurachSoldatov2016MFAT}
\emph{Mikhailets V. A., Murach A. A., Soldatov V.} A criterion for continuity in a parameter of solutions to generic boundary-value problems for higher-order differential systems~// Methods Funct. Anal. Topology.~-- 2016.~-- Vol.~22, No.~4.~-- P. 375\,--\,386.

\bibitem{Gnyp2016}
\emph{Gnip (Gnyp) E. V.} Continuity with respect to the parameter of the solutions of one-dimensional boundary value problems in Slobodetskii spaces~// Ukrainian Math.~J.~-- 2016.~-- \textbf{68}, №~6.~-- P. 849\,--\,861.

\bibitem{GnypMikhailetsMurach2017}
\emph{Hnyp (Gnyp) E., Mikhailets V., Murach A.} Parameter-dependent one-dimensional boundary-value problems in Sobolev spaces~// Electron.~J. Differential Equations.~-- 2017.~-- №~81.~-- P. 1\,--\,13.


\bibitem{MasliukPelekhataSoldatov2020MFAT}
\emph{Masliuk H., Pelekhata O., Soldatov V.} Approximation properties of multipoint boundary-value problems~// Methods Funct. Anal. Topology.~-- 2020.~-- \textbf{26}, №~2.~-- P. 119\,--\,125.

\bibitem{DanfordShvarts1958}
\emph{Dunford N., Schwartz J. T.} Linear Operators. Part I: General Theory.~--  New York: Interscience, 1958.~-- xiv+858~p.

\bibitem{ReedSimon}
\emph{Reed M., Simon B.} Methods of Modern Mathematical Physics. I. Functional analysis.~-- New York: Academic Press, 1980.~-- xv+400~p.

\bibitem{RieszSz-Nagy90}
\emph{Riesz F., Sz-Nagy B.} Functional Analysis.~-- New York: Dover Publ. Inc., 1990.~-- xii+504~p.

\end{thebibliography}
\end{document}